%% file: nonexist.tex
\documentclass[11pt]{article}
\usepackage[francais,english]{babel}
\usepackage[latin1]{inputenc}

\usepackage{amsmath}
\usepackage{amsfonts}
\usepackage{amssymb}
\usepackage{amsthm}
\usepackage{graphics}
\usepackage{amscd}
\usepackage{graphicx,epsfig}

\DeclareMathOperator{\Div}{div}

\renewcommand{\epsilon}{\varepsilon}

\newcommand{\boB}{\mathcal{B}}

\newcommand{\Ome}{\Omega}

\newcommand{\R}{\mathbb{R}}

\newcommand{\C}{\mathbb{C}}

\newcommand{\der}[2]{\dfrac{\partial #1}{\partial #2}}
\newcommand{\dd}{\mathrm{d}}

\newtheorem{thm}{Th\'eor\`eme}

\newtheorem{prop}{Proposition}
\newtheorem{lem}{Lemme}

\renewcommand{\phi}{\varphi}
\newcommand{\dis}{\displaystyle}

\title{Quelques r\'esultats de non-existence pour l'\'equation des surfaces
  minimales}
\author{Laurent Mazet}
\date{}

\begin{document}

\maketitle

\selectlanguage{french}
\begin{abstract}
Dans cet article, nous d\'emontrons qu'il n'existe pas de solution $u$ de
l'\'equation des surfaces minimales sur un domaine asymptotiquement \'egal
au secteur angulaire valant $+\infty$ sur un de ses bords et $-\infty$ sur
l'autre. 
\end{abstract}

\selectlanguage{english}

\begin{abstract}
In this article, we prove that there does not exist a solution $u$ of the
minimal surfaces equation on a domain, which is asymptoticaly an angular
sector, and taking the value $+\infty$ on one side and $-\infty$ on the
other. 
\end{abstract}

\selectlanguage{french}
\noindent 2000 \emph{Mathematics Subject Classification.} 53A10.

\section*{Introduction}
\input nonexist0.tex
\section{Quelques cas particuliers}{\label{S1}}
\input nonexist1.tex
\section{Le r\'esultat g\'en\'eral}
\input nonexist2.tex
\section{Une g\'en\'eralisation au cas $\alpha\ge \pi$}
\input nonexist3.tex

\bigskip

\noindent Laurent Mazet

\noindent Laboratoire Emile Picard (UMR 5580), Universit\'e Paul Sabatier,

\noindent 118, Route de Narbonne, 31062 Toulouse, France.

\noindent E-mail: mazet@picard.ups-tlse.fr

\end{document}

%% file: nonexist0.tex
Une fonction $u$ d\'efinie sur un domaine $\Ome$ de $\R^2$ satisfait
l'\'equation des surfaces minimales si: 

\begin{equation*}
\Div\left(\frac{\nabla u}{\sqrt{1+|\nabla u|^2}} \right)=0 \tag{ESM}
\label{MSE}
\end{equation*}
Cette \'equation aux d\'eriv\'ees partielles traduit le fait que le graphe de
la fonction $u$ est une surface minimale de $\R^3$. A cette \'equation, on
associe le probl\`eme de Dirichlet: il s'ag\^\i t de trouver une solution $u$
de 
l'\'equation des surfaces minimales sur un domaine $\Ome$ en imposant la
valeur de $u$ sur le bord de $\Ome$. 

Plusieurs solutions de l'\'equation \eqref{MSE} sont connues; ainsi sur la
demi-bande 
$\R_+ \times [-\frac{\pi}{2},\frac{\pi}{2}]$, la fonction $h(x,y)=x\tan y$ est
solution de \eqref{MSE}. Le graphe de $h$ est un morceau d'h\'elico\" \i
de. On constate que si $(x_n,y_n)$ converge vers un point de $\R_+^*\times
\{\frac{\pi} {2}\}$ la suite $(h(x_n,y_n))$ converge vers $+\infty$ et si
$(x_n,y_n)$ converge vers un point de $\R_+^*\times \{-\frac{\pi} {2}\}$, on a
alors $\lim h(x_n,y_n)=-\infty$. $h$ se pr\'esente donc comme une solution du
probl\`eme de Dirichlet sur la demi-bande avec les valeurs $0$ sur $\{0\}\times
]-\frac{\pi}{2},\frac{\pi}{2}[$, $+\infty$ sur $\R_+^*\times \{\frac{\pi}
{2}\}$ et $-\infty$ sur $\R_+^*\times \{-\frac{\pi} {2}\}$.

Maintenant la question que l'on se pose est de savoir ce qu'il se passe si les
deux demi-droites appartenant au bord de la demi-bande ne sont plus
parall\`eles mais forment un angle $2\alpha>0$. Plus pr\'ecis\'ement, si
$\Ome$ est le domaine $\{(x,y)\in [1,+\infty[\times \R|
-\tan\alpha<\frac{y}{x}<\tan \alpha \}$, existe-t'il une solution $u$ de
\eqref{MSE} sur $\Ome$ telle que $u$ prenne la valeur $0$ sur $\Ome\cap
\{x=1\}$, tende 
vers $+\infty$ sur $\Ome\cap \{x\sin\alpha-y\cos\alpha=0\}$ et tende vers
$-\infty$ sur $\Ome\cap \{x\sin\alpha+y\cos\alpha=0\}$?

On connait de nombreux r\'esultats concernant le problème de Dirichlet sur le
secteur ainsi, dans \cite{RSE}, H.~Rosenberg et R.~Sa-Earp construisent des
solutions 
pour des donn\'ees continues sur le bord d'un domaine convexe inclus dans un
secteur angulaire. Dans \cite{Ni}, J.C.C.~Nitsche se pose la question d'un
principe du maximum pour le secteur, un tel principe est d\'emontr\'e dans
\cite{RSE} ce qui implique, d'après le r\'esultat de P.~Collin et R. Krust
dans \cite{CK}, l'unicit\'e des 
solutions ayant des donn\'ees born\'ees sur le bord. On a
d'autres r\'esultats d'unicit\'e, ainsi C.-C.~Lee, dans \cite{Le},
d\'emontre l'unicit\'e en imposant des conditions sur les d\'eriv\'ees au
bord. 

La r\'eponse à la question pos\'ee plus haut est apport\'ee par le Th\'eorème
\ref{T1} et cette r\'eponse est n\'egative: une solution de \eqref{MSE}
prenant la valeur $+\infty$ sur un cot\'e d'un secteur et $-\infty$ sur
l'autre ne peut exister. 

Le Th\'eorème \ref{T1} est le principal r\'esultat de cet article. Il
concerne les domaines $\Ome$ tels que, hors d'un disque de $\R^2$, $\Ome=\{r\cos
\theta, r\sin \theta)\in\R^2 |\ r>r_0,\ -\alpha<\theta<\alpha\}$ avec
$0<\alpha<\pi$. Le Th\'eorème \ref{T1} affirme alors qu'il n'existe pas de
solution $u$ de \eqref{MSE} sur $\Ome$ telle que $u$ tende vers $+\infty$ sur
l'un des cot\'es de ce secteur angulaire et $-\infty$ sur l'autre.

\medskip

Dans une première partie, nous donnons quelques r\'esultats
pr\'eliminaires. Parmi ceux-ci le plus important est la Proposition \ref{P3}
qui donne le 
r\'esutat du Th\'eorème \ref{T1} dans le cas de petits angles $\alpha$.

La deuxième partie est consacr\'ee à la preuve du Th\'eorème
\ref{T1}. L'id\'ee principale de la preuve est que l'existence de fonctions
$u$ contredisant le r\'esultat pour de grands angles implique l'existence de
fonctions $v$ contredisant le r\'esultat pour de petits angles. On entrerait
ainsi en contradiction avec la Proposition \ref{P3}.

La troisième partie donne une g\'en\'eralisation de ce r\'esultat dans le cas
où l'angle $\alpha$ est sup\'erieur à $\pi$.

\medskip

Dans la suite, si $u$ est une fonction d\'efinie sur un domaine $\Ome$ on
notera $W=\sqrt{1+|\nabla u|^2}$.

%%% Local Variables: 
%%% mode: latex
%%% TeX-master: t
%%% TeX-master: t
%%% End: 

%% file: nonexist1.tex
\subsection{Pr\'eliminaires}
\begin{prop}{\label{P1}}
Soit $\Omega$ un domaine de $\R^2$ tel qu'une composante connexe
de son bord soit une droite $L$; on suppose de plus qu'il existe une
partie $\Delta$ de $\R^2$ telle que $\Omega \backslash \Delta$ soit
isom\'etrique \`a la bande $]-\epsilon,0[\times \R$ ($\epsilon>0$)
o\`u $L$ est la droite $\{0\}\times\R$. Alors il n'existe pas de
solution $u$ de l'\'equation des surfaces minimales qui
prend la valeur $+\infty$ sur $L$ (de m\^eme pour $-\infty$).
\end{prop}

\begin{proof}
Supposons, par l'absurde, que $u$ soit une solution de l'\'equation
des surfaces minimales sur $]-\epsilon,0[\times \R$ valant $+\infty$
sur $L$. Alors le graphe de $u$ contredit le principe du demi-espace,
plus pr\'ecis\'ement, la d\'emonstration du Th\'eor\`eme $1$ dans \cite{HM} 
contredit l'existence du graphe de $u$.  
\end{proof}

Par la suite nous allons avoir besoin d'estim\'ees sur les
d\'eriv\'ees de solution de l'\'equation des surfaces minimales. On a
alors le r\'esultat suivant.

\begin{lem}{\label{jenkins}}
Soit $\Ome$ un domaine convexe du plan $\R^2$ et $u$ une solution de
l'\'equation des surfaces minimales sur $\Ome$. On consid\`ere $p$ un
point de $\Ome$ et $d$ la distance du point $p$ au bord de $\Ome$. On
note $q\in\partial \Ome$ un point qui r\'ealise cette distance, $n$
le vecteur unitaire $\dis\frac{\overrightarrow{pq}}{d}$ et $n'$ un vecteur
unitaire normal \`a $n$. On note
maintenant $\Sigma$ le graphe de $u$ et $P$ le point de $\Sigma$ au
dessus de $p$. Soit $r$ la distance le long de la surface $\Sigma$ du
point $P$ au bord de $\Sigma$. Alors si le rapport $\dis\frac{d}{r}$ est
inf\'erieur \`a $\dis\frac{1}{8}$, on a au point $p$:
\begin{gather*}
\frac{|n\cdot \nabla u(p)|}{W}\ge 1-4\frac{d^2}{r^2}\\
\frac{|n'\cdot \nabla u(p)|}{W} \le 2\sqrt{2}\frac{d}{r}
\end{gather*}
\end{lem}

Il s'agit du Lemme $1$ dans \cite{JS}.

\subsection{Le cas du secteur angulaire}

Pour $0<\alpha<\pi$, on consid\`ere le domaine de $\R^2$:
$$D(\alpha)=\{(r\cos\theta,
r\sin\theta)\}_{r\in\R_+,\theta\in[-\alpha,\alpha]}$$
Pour tout $\theta$, on note $L(\theta)=\{(r\cos\theta,
r\sin\theta)\}_{r\in\R_+^*}$; on a alors la proposition suivante:

\begin{prop}{\label{P2}}
Soit $\alpha\in ]0,\pi[$. Il n'existe pas de solution $u$ de
    l'\'equation des surfaces minimales sur $D(\alpha)$ telle que $u$
    prenne la valeur $+\infty$ (resp. $-\infty$) sur $L(\alpha)$
    (resp. $L(-\alpha)$). 
\end{prop}  

\begin{proof}
Supposons qu'une telle solution $u$ existe. Alors d'apr\`es le
Th\'eor\`eme $3$ et la Remarque $1$ dans \cite{Ma1}, on peut prolonger le
graphe de la fonction $u$ 
par sym\'etrie par rapport \`a la droite verticale de $\R^3$:
$L=\{x=0\}\cap\{y=0\}$. La surface $\Sigma$ ainsi obtenue est une surface
minimale compl\`ete simplement connexe ayant un nombre fini de points
de branchement le long de la droite verticale de sym\'etrie $L$. Comme
$\Sigma$ est un graphe au dessus de $D(\alpha)$, les points de
$\Sigma$ ayant une normale horizontale sont les points de $L$

D'apr\`es le Lemme $4$ dans \cite{Ma2}, la courbe $r\mapsto (r\cos\theta,
r\sin\theta, u(r\cos\theta,r\sin\theta))$ appartenant \`a $\Sigma$ a pour 
extr\'emit\'e lorsque $r\rightarrow 0$ un point de $L$ o\`u la normale
vaut $(\sin\theta,-\cos\theta,0)$ ou $(-\sin\theta,\cos\theta,0)$. De
plus tout point de $L$ est l'extr\'emit\'e d'une telle courbe. Ceci
implique que pour $\theta\in]\alpha-\pi,\pi-\alpha[$ l'un des
    vecteurs suivant $(\sin\theta,-\cos\theta,0)$ ou
    $(-\sin\theta,\cos\theta,0)$ n'est jamais la normale \`a $\Sigma$,
    en effet ces deux vecteurs ne peuvent \^etre la normale \`a $\Sigma$
    qu'en l'extr\'emit\'e de $r\mapsto (r\cos\theta, r\sin\theta,
u(r\cos\theta,r\sin\theta))$ ce qui ne fait qu'un seul point. On sait
ainsi que la normale \`a $\Sigma$ omet un nombre infini de points de la
sph\`ere et m\^eme un segment dans la sph\`ere par continuit\'e de la
normale. Or ceci est impossible pour une surface minimale compl\`ete
simplement connexe ayant un nombre fini de points de branchement (Th\'eor\`eme 
$8.2$ dans \cite{Os}). 
\end{proof}

\subsection{Le cas des petits angles}

Soit $\alpha\in]0,\frac{\pi}{2}[$ et $a$ un nombre r\'eel positif. On
note alors $D_a(\alpha)=D(\alpha)\cap \{x\ge a\}$. Par abus de
notation, on continue \`a noter $L(\theta)$ l'intersection de
$L(\theta)$ avec $D_a(\alpha)$. 

\begin{prop}{\label{P3}}
Il existe $0<\alpha_0<\frac{\pi}{2}$ tel que si $0<\alpha<\alpha_0$ et
$a\in \R_+$, il n'existe pas de solution $u$ de l'\'equation des surfaces
minimales sur $D_a(\alpha)$ telle que $u$ prenne la valeur $+\infty$
(resp. $-\infty$) sur $L(\alpha)$ (resp. $L(-\alpha)$). 
\end{prop}

\begin{proof}
Supposons qu'une telle solution $u$ existe. On consid\`ere $p$ un
point de $D_a(\alpha)$ de coordonn\'ees $(x,y)$. Soit $x_0>a$ tel que
$x_0\sin\alpha=x_0-a$. Si le point $p$ v\'erifie $x\ge x_0$, le point
$q$ du bord de $D_a(\alpha)$ qui r\'ealise la distance de $p$ au bord
de $D_a(\alpha)$ n'est pas situ\'e sur $\{x=a\} $ mais sur
$L(\alpha)$ ou $L(-\alpha)$ et on a $|pq|\le x\sin\alpha$. Le long du
graphe de $u$ la distance de $P$, le point de coordonn\'ees
$(x,y,u(x,y))$, au bord du graphe est minor\'ee par $x-a$. Le rapport
entre ces deux longueurs est donc major\'ee par 
$$\frac{x\sin\alpha}{x-a}=\sin\alpha \frac{1}{1-\frac{a}{x}}$$
Posons $x_1=\max\{x_0,10a\}$, alors si $p$ v\'erifie $x\ge x_1$ le
rapport entre ces deux longueurs est major\'ee par $\frac{10}{9}\sin
\alpha$. Notons $\alpha_1\in[0,\frac{\pi}{2}]$ le nombre v\'erifiant
$\frac{10}{9} \sin\alpha_1=\frac{1}{8}$. Donc si $\alpha<\alpha_1$ on
peut appliquer le Lemme \ref{jenkins} au point $p$. Le vecteur
$\dis\frac{\overrightarrow{pq}}{|pq|}$ est $(-\sin\alpha, \cos\alpha)$ ou
$(-\sin \alpha, -\cos\alpha)$. On en d\'eduit alors, par le Lemme
\ref{jenkins}  que
\begin{equation*}
\begin{split}
\frac{|u_y|}{W}&\ge(1-4\left(\frac{10}{9}\right)^2\sin^2 \alpha)\cos
\alpha- 2\sqrt{2}\frac{10}{9}\sin^2 \alpha\\
&\ge \cos\alpha-\sin^2 \alpha \left(4\left(\frac{10}{9}\right)^2\cos
\alpha+ 2\sqrt{2} \frac{10}{9}\right) 
\end{split}
\end{equation*}
Il existe donc $\alpha_0<\alpha_1$ tel que si $\alpha\le \alpha_0$ on
a $|u_y(p)|\neq 0$, pour tout point $p=(x,y)$ avec $x\ge x_1$, et donc $u_y>0$
(On remarque que $\alpha_0$ ne d\'epend pas de $a$).

On se place maintenant avec $\alpha\le \alpha_0$ et sur
$D_a(\alpha)\cap \{x\ge x_1\}$ on pose $\Phi(x,y)=(x,u(x,y))$. On
d\'efinit ainsi une application \`a valeur dans $[x_1,+\infty[\times
\R$. La diff\'erentielle de $\Phi$ est 
$$\dd \Phi(x,y)=
\begin{pmatrix}
1& u_x\\
0& u_y
\end{pmatrix}$$ 
dont le jacobien est non nul puisque sa valeur est $u_y\neq 0$. $\Phi$
est donc un diff\'eomorphisme local; en fait, il s'agit d'un
diff\'eomorphisme global de $D_a(\alpha)\cap \{x\ge x_1\}$ sur
$[x_1,+\infty[\times \R$ car $u$ croit strictement de $-\infty$ \`a
$+\infty$ le long des courbes $x=\textrm{cste}$. 

On voit ainsi que le graphe de $u$ au dessus de $D_a(\alpha)\cap
\{x\ge x_1\}$ peut \^etre vu comme le graphe d'une fonction $v$
solution de \eqref{MSE} au dessus du demi-plan vertical $\{y=0,x\ge
x_1\}$. La fonction $v$ est alors born\'ee sur $\{x=x_1\}$; donc,
d'apr\`es le r\'esultat de P.~Collin et R.~Krust (Th\'eorème $3.4$ dans
\cite{CK}), le graphe est  
asymptote \`a un plan d'equation $y=cx+d$ et de plus ils existent
$d_1$ et $d_2$ tel que l'on ait $cx+d_1 \le
v(x,z)\le cx+d_2$ ce qui est imposible \`a cause de la forme de
$D_a(\alpha)$. 
\end{proof}

%%% Local Variables: 
%%% mode: latex
%%% TeX-master: t
%%% End: 

%% file: nonexist2.tex
Nous avons maintenant un r\'esultat similaire \`a celui de la
Proposition \ref{P3} mais sans limite sur l'angle $\alpha$.

\begin{thm}{\label{T1}}
Soit $\Ome$ un domaine de $\R^2$ tel qu'il existe un compact $K$ de
$\R^2$ et $\alpha\in]0,\pi[$ v\'erifiant $\Ome \backslash K=D(\alpha)
\backslash D(0,r)$ o\`u $D(0,r)$ d\'esigne le disque de centre
l'origine et de rayon $r$. Alors, il n'existe pas de solution $u$ de
l'\'equation des surfaces minimales sur $\Ome$ telle que $u$ prenne la
valeur $+\infty$ (resp. $-\infty$) sur $L(\alpha)$ (resp.
$L(-\alpha)$).  
\end{thm}

\begin{proof}
Supposons qu'une telle solution existe, nous avons donc une fonction
$u$ sur $D(\alpha) \backslash D(0,1)$ (on suppose $r=1$) qui prend la
valeur $+\infty$ (resp. $-\infty$) sur $L(\alpha)$
(resp. $L(-\alpha)$). L'id\'ee de la preuve est de se servir de $u$
pour construire une fonction qui contredirait la Proposition \ref{P3}.

Soit $\beta>0$ tel que $\beta$ est inf\'erieur \`a $\alpha$ et
$\frac{\pi}{4}$. On note alors $A_r$, $B_r$, $C_r$ et $D_r$ les points
de $D(\alpha)$ de coordonn\'ees polaires respectives $(r,\alpha)$,
$(r,\alpha-\beta)$, $(r,-\alpha+\beta)$ et $(r,-\alpha)$. On a alors
le r\'esultat suivant

\begin{lem}{\label{L1}}
Il existe $\beta\le \min\{\alpha,\frac{\pi}{4}, 2\alpha_0\}$
($\alpha_0$ est l'angle donn\'e par la Proposition \ref{P3}) et $r_0>1$
tels que la fonction $u$ soit croissante
\begin{itemize}
\item le long des deux arcs de cercles de
centre l'origine joignant $B_r$ \`a $A_r$ et $D_r$ \`a $C_r$ et
\item le long des deux segments $[B_r,A_r]$ et $[D_r,C_r]$
\end{itemize}
pour $r>r_0$.
\end{lem}

\begin{proof}
Par sym\'etrie du probl\`eme, on ne consid\`erera que le segment
$[B_r,A_r]$ et l'arc de cercle joignant $B_r$ \`a $A_r$. Soit $p$ un
point de coordonn\'ees polaires $(\rho,\theta)$ avec
$\theta\in[\alpha-\beta,\alpha]$ et $\rho\ge \frac{2}{2-\sqrt{2}}$; alors
le point $q$ du bord de $D(\alpha)\backslash D(0,1)$ qui r\'ealise la
distance de $p$ au bord appartient \`a $L(\alpha)$. La distance $|pq|$
est donc major\'ee par $\rho\sin\beta$, la distance le long du graphe
de $u$ de $P$ le point du graphe situ\'e au dessus de $p$ au bord est
minor\'ee par $\rho-1$. Le rapport entre ces deux longueurs est donc
major\'e par $\sin\beta \frac{\rho}{\rho-1}$. Ainsi, pour $\rho\ge 10$, le
rapport est major\'e par $\frac{10}{9}\sin\beta$; soit $\beta_1$ tel
que $\frac{10}{9}\sin\beta_1=\frac{1}{8}$ alors pour $\beta<\beta_1$
le rapport est inf\'erieur \`a $\frac{1}{8}$.

On suppose $\beta<\beta_1$. On d\'efinit $r_0$ tel que
$r_0\cos\frac{\pi}{8}=10$. Soit $r\ge r_0$ on suppose maintenant que
$p$ appartienne au segment $[B_r,A_r]$ on a alors $\rho\ge 10$ par
d\'efinition de $r_0$. Le
vecteur $\frac{\overrightarrow{pq}}{|pq|}$ est
$(-\sin\alpha,\cos\alpha)$ et le vecteur unitaire directeur de
$[B_r,A_r]$ est $n=(-\sin(\alpha-\frac{\beta}{2}),
\cos(\alpha-\frac{\beta}{2}))$. D'apr\`es le Lemme \ref{jenkins}, on a donc
\begin{equation*}
\frac{|\nabla u\cdot n|}{W}(p)\ge (1-4(\frac{10}{9}\sin\beta)^2) 
(-\sin\alpha,\cos\alpha)\cdot n -2\sqrt{2}\frac{10}{9}\sin\beta
\end{equation*}
On constate que losque $\beta$ tend vers $0$ le minorant tend vers
$1$. Il existe donc $0<\beta_2<\beta_1$ tel que pour tout
$\beta<\beta_2$ : $\frac{|\nabla u\cdot n|}{W}(p)>0$.

Si maintenant $p$ est un point de l'arc de cercle joignant $B_r$ \`a
$A_r$ le vecteur unitaire tangent \`a l'arc de cercle en $p$ est
$n'(\theta)=(-\sin\theta, \cos\theta)$. D'apr\`es le Lemme
\ref{jenkins}, on a donc
\begin{equation*}
\begin{split}
\frac{|\nabla u\cdot n'(\theta)|}{W}(p)&\ge (1-4(\frac{10}{9}\sin\beta)^2) 
(-\sin\alpha,\cos\alpha)\cdot n'(\theta) -2\sqrt{2}\frac{10}{9}\sin\beta\\ 
&\ge (1-4(\frac{10}{9}\sin\beta)^2) (-\sin\alpha,\cos\alpha) \cdot
n'(\beta)- 2\sqrt{2}\frac{10}{9}\sin\beta 
\end{split}
\end{equation*}
De m\^eme que pr\'ec\'edemment, le minorant tend vers $1$ lorsque
$\beta$ tend vers $0$. Il existe donc $0<\beta_3<\beta_2$ tel que pour
tout $\beta<\beta_3$ : $\frac{|\nabla u\cdot n'(\theta)|}{W}(p)>0$.
Alors tout $\beta<\min\{ \beta_3,2\alpha_0\}$ r\'epond au lemme.
\end{proof}

On fixe maintenant $\beta$ et $r_0$ tels que le Lemme \ref{L1} soit
satisfait. On a alors le r\'esultat suivant
\begin{lem}{\label{L2}}
Il existe $r_1\ge r_0$ tel que si $p\in D(\alpha)$ est le point de
coordonn\'ees polaires $(r,\theta)$ avec $r\ge r_1$ on a
$\der{u}{\theta}(p)>0$. 
\end{lem}

\begin{proof}
Tout d'abord le Lemme \ref{L1} nous dit que la d\'eriv\'ee est
positive si $\theta\in ]-\alpha,-\alpha+\beta]\cup [\alpha-\beta,
\alpha[$. Maintenant si le lemme n'est pas vrai il existe une suite de
point $p_n$ de coordonn\'ees polaires $(r_n,\theta_n)$ telle que
$r_n\rightarrow+\infty$ et $\der{u}{\theta}(p_n)\le 0$. On a
$\theta_n\in 
[-\alpha+\beta,\alpha-\beta]$ et, quitte \`a extraire, on peut supposer
que $\theta_n\rightarrow \theta_\infty\in
[-\alpha+\beta,\alpha-\beta]$. On d\'efinit alors $\Ome_n$ le 
domaine de $\R^2$: $\Ome_n=\{(r\cos\theta,r\sin\theta)| \
r\ge{r_n}^{-1},\ \theta \in]-\alpha+\theta_\infty-\theta_n,
\alpha+\theta_\infty-\theta_n[\} $. Sur $\Ome_n$ on d\'efinit, en
coordonn\'ees polaires, la fonction $v_n$: $v_n(r,\theta)= {r_n}^{-1}
u(r_nr,\theta+\theta_\infty-\theta_n)$. $v_n$ est une solution de
l'\'equation des surfaces minimales qui vaut $+\infty$ sur
$L(\alpha+\theta_\infty-\theta_n)$ et $-\infty$ sur $L(-\alpha+
\theta_\infty-\theta_n)$. L'hypoth\`ese sur $u$ au point $p_n$ se
traduit par $\der{v_n}{\theta}(1,\theta_\infty)\le 0$. 

On va maintenant \'etudier la limite de la suite $(v_n)$; pour cela on va
\'etudier les lignes de divergence de $(v_n)$ (voir \cite{Ma1} et
\cite{Ma2}). Le domaine
limite est $D(\alpha)$. Comme $v_n$ prend des valeurs infinies sur
$L(\alpha+\theta_\infty-\theta_n)$ et $L(-\alpha+ \theta_\infty-
\theta_n)$, les seules lignes de divergence possibles sont les droites
incluses dans $D(\alpha)$ et les demi-droites ayant l'origine pour
extr\'emit\'e; on utilise ici le Lemme A.$1$ de \cite{Ma2}.

Consid\'erons $L(\gamma)\subset D(\alpha)$ et supposons que $L(\gamma)$
soit une ligne de divergence. On consid\`ere $\Psi_n$ la fonction
conjugu\'ee de $v_n$ normalis\'ee (en coordonn\'ees polaires) par
$\Psi_n ({r_n}^{-1},\theta_\infty)=0$. Alors en
choisissant la bonne sous-suite, on a $\Psi_n\rightarrow \Psi$ avec
$\Psi$ qui v\'erifie: $\Psi=0$ en l'origine,
$\Psi(\rho,\alpha)=-\rho$, $\Psi(\rho,-\alpha)=-\rho$ et
$\Psi(\rho,\gamma)=\pm \rho$ suivant la valeur de la normale
limite. Mais si $\gamma\ge 0$, comme $\Psi$ est $1$-lipschitzienne on a
$|\Psi(\rho,\alpha)- \Psi(\rho,\gamma)|<2\rho$; donc
$\Psi(\rho,\gamma)=-\rho$ et la normale limite le long de $L(\gamma)$
est $(-\sin\gamma,\cos\gamma)$; de m\^eme si $\gamma\le 0$. Ainsi pour
$L(\gamma)$ il n'y a qu'une possibilit\'e pour la normale limite.  

Supposons que la droite $L\subset D(\alpha)$ soit une
ligne de divergence; on a alors $N_{n'}$ qui converge vers une normale
\`a $L$ le long de $L$. Soit $L'$ la parall\`ele \`a $L$ passant par
l'origine. On s'int\'eresse maintenant \`a la suite $(v_{n'})$. Si une
ligne de divergence pour $(v_{n'})$ existe dans la bande comprise
entre $L$ et $L'$ ce 
doit \^etre une droite parall\`ele \`a $L$. Si un point de la bande
comprise entre $L$ et $L'$ appartient au domaine de convergence de
$(v_{n'})$, la composante connexe de $\boB(v_{n'})$ qui contient ce
point 
v\'erifie alors les hypoth\`eses de la Proposition \ref{P1} et la limite
de la suite $v_{n'}$ sur cette composante prend une valeur infinie sur
une droite parall\`ele \`a $L$; ceci est impossible par la Proposition
\ref{P1}. Ainsi en tout point de la bande comprise entre $L$ et $L'$
passe une ligne de divergence; de plus en tout point de cette bande
$N_{n'}$ converge vers le m\^eme vecteur unitaire que le long de
$L$, ceci est entre autre vrai le long de $L'$ sauf pour
l'origine. Ainsi il existe $\gamma$ tel que $L(\gamma)$ et $L(\gamma-\pi)$
soient deux lignes de divergence avec la m\^eme normale limite ce qui
est impossible d'apr\`es le r\'esultat du paragraphe pr\'ec\'edent. Donc $L$ ne
peut \^etre une 
ligne de divergence. Et les $L(\gamma)$ sont les seules lignes de
divergence possibles.

Si le point de coordonn\'ees polaires $(1,\theta_\infty)$ est un point
d'une ligne de divergence la normale limite \`a ce point est $(-\sin
\theta_\infty, \cos \theta_\infty)$ ce qui contredit que
$\der{v_n}{\theta}(1,\theta_\infty)\le 0$. Ce point appartient donc
au domaine de convergence, or la composante connexe du domaine de
convergence le contenant est un secteur angulaire compris
entre $L(\gamma_1)$ et $L(\gamma_2)$ (avec $-\alpha\le \gamma_1<\gamma_2
\le \alpha$) et la limite $v$ de la suite $(v_n)$ sur cette
composante vaut $+\infty$ sur $L(\gamma_2)$ et $-\infty$ sur
$L(\gamma_1)$ d'apr\`es les normales limites de long de ces deux lignes
de divergence (Lemme A.2 de \cite{Ma2}). Une telle solution $v$ est impossible
par la 
Proposition \ref{P2}. On vient ainsi d'\'etablir une contradiction,
le lemme est donc d\'emontr\'e.
\end{proof}

Nous pouvons maintenant finir la preuve du th\'eor\`eme. On
consid\`ere le domaine $D_{r_1\sin\frac{\beta}{2}} (\frac{\beta}{2})$  
(avec $\beta$ donn\'e par Lemme \ref{L1} et $r_1$ donn\'e par le Lemme 
\ref{L2}) et on note $E_r$ le point de coordonn\'ees polaires $(r,
\frac{\beta}{2})$ et $F_r$ le point de coordonn\'ees polaires $(r,
-\frac{\beta}{2})$. Soit $n$ un entier positif, les polygones
 $A_{r_1}B_{r_1}B_{r_1+n}A_{r_1+n}$, $C_{r_1}D_{r_1}
 D_{r_1+n}C_{r_1+n}$ et 
 $E_{r_1}F_{r_1} F_{r_1+ n}E_{r_1+n}$ se d\'eduisent alors les uns des
 autres  par des d\'eplacements et peuvent \^etre ainsi identifi\'es. En
 utilisant ces identification  on construit sur 
 $E_{r_1}F_{r_1}  F_{r_1+n}E_{r_1+n}$ la solution $w_n$ de l'\'equation des
 surfaces minimales valant $u_{|[A_{r_1},B_{r_1}]}$ sur
 $[E_{r_1},F_{r_1}]$, $-\infty$ sur $[F_{r_1}, F_{r_1+ n}]$,
 $u_{|[A_{r_1+n},B_{r_1+n}]}$ sur $[E_{r_1+n},F_{r_1+n}]$ et $+\infty$ 
 sur $[E_{r_1}, E_{r_1+ n}]$ ($w_n$ existe d'apr\`es le Th\'eor\`eme $2$ de \cite{JS}). D'apr\`es le principe du maximum et le
 Lemme \ref{L2}, on a alors:
$$u_{|C_{r_1}D_{r_1} D_{r_1+n}C_{r_1+n}} \le w_n \le
u_{|A_{r_1}B_{r_1}B_{r_1+n}A_{r_1+n}}$$ 

Ainsi $(w_n)$ est uniform\'ement born\'ee sur tout compact de
$D_{r_1\sin\frac{\beta}{2}} (\frac{\beta}{2})$. Donc une sous suite $(w_{n'})$
converge sur $D_{r_1\sin\frac{\beta}{2}} 
(\frac{\beta}{2})$ vers une solution $w$ de \eqref{MSE} prenant la valeur
$+\infty$ sur $L(\frac{\beta}{2})$ et $-\infty$ sur
$L(-\frac{\beta}{2})$; ainsi $w$ contredit la Proposition \ref{P3} et le
th\'eor\`eme est d\'emontr\'e.   
\end{proof}

%%% Local Variables: 
%%% mode: latex
%%% TeX-master: t
%%% End: 

%% file: nonexist3.tex
On souhaite pouvoir donn\'e un sens au probl\`eme d'existence lorsque
$\alpha$ est plus grand que $\pi$ et voir si le Th\'eor\`eme \ref{T1}
reste vrai dans ce cas.

Pour cela posons $D(\alpha_1, \alpha_2)=\{(r,\theta)|\ r\ge 0,
\alpha_1\le 
\theta\le \alpha_2\}$ muni de la m\'etrique polaire $\dd s^2=\dd
r^2+r^2\dd\theta^2$ (tout les points de la forme $(0,\theta)$ sont
identifi\'es, on appelle $0$ ce nouveau point: c'est le sommet de
$D(\alpha_1,\alpha_2)$). Soit $\phi$ l'application d\'efinie sur
$D(\alpha_1,\alpha_2)$ par 
$(r,\theta)\mapsto (r\cos\theta, r\sin\theta)$, $\phi$ est alors
une isom\'etrie locale de $D(\alpha_1,\alpha_2)$ sur $\R^2$. Le couple 
$(D(\alpha_1,\alpha_2),\phi)$ est un multi-domaine (voir \cite{Ma1}). Pour
$\alpha>0$, 
on note $D(\alpha)=D(-\alpha,\alpha)$; on constate que, pour
$\alpha<\pi$, cette d\'efinition est \'equivalente \`a celle de la
section \ref{S1}. Comme pr\'ec\'edemment, on note
$L(\theta)=\{(r,\theta)|\ r>0\}$.

Pour pouvoir g\'en\'eraliser le Th\'eor\`eme \ref{T1}, nous allons
devoir rajouter des hypoth\`eses.

\subsection{Le cas du secteur angulaire}

Dans cette partie, nous g\'en\'eralisons la Proposition \ref{P2} mais
pour cela nous devons rajouter une hypoth\`ese sur $\Psi_u$. 
\begin{prop}{\label{P4}}
Soit $\alpha\ge\pi$. Il n'existe pas de solution $u$ de l'\'equation
des 
surfaces minimales sur $D(\alpha)$ telle que $u$ prenne la valeur
$+\infty$ (resp. $-\infty$) sur $L(\alpha)$ (resp. $L(-\alpha)$) et telle que
$\Psi_u\le 0$ ($\Psi_u$ est normalis\'ee par $\Psi_u(0)=0$) 
\end{prop}

\begin{proof}
Supposons qu'une telle solution $u$ existe. De m\^eme que dans la
preuve de la Proposition \ref{P2}, on peut prolonger le graphe de $u$
par 
sym\'etrie par rapport \`a l'axe verticale $L=\{x=0,y=0\}$. La
surface $\Sigma$ ainsi obtenue est une surface minimale compl\`ete de
$\R^3$, elle est
simplement connexe et r\'eguli\`ere car $\Psi_u\le 0$.
% et de courbure totale finie.
On sait de plus que la normale \`a $\Sigma$ est horizontale uniquement le long
de $L$. Par le Lemme 4 dans \cite{Ma2}, le vecteur horizontal
$(\sin\theta,-\cos\theta,0)$ ne peut \^etre la normale \`a $\Sigma$ qu'en
l'extremit\'e de la courbe incluse dans $\Sigma$: $r\mapsto
(r\cos(\theta+k\pi), r\sin(\theta+k\pi), u(r,\theta+k\pi))$ avec
$-\alpha<\theta+k\pi<\alpha$. Ainsi tout vecteur unitaire horizontal est
normal \`a $\Sigma$ qu'en un nombre fini de point. Ces vecteurs s'identifient
\`a l'\'equateur de la sph\`ere unit\'e, or l'\'equateur est un
ensemble de capacit\'e logarithmique non nulle, $\Sigma$ est donc de courbure
totale finie d'apr\`es le Corollaire du Th\'eor\`eme 2.3 dans
\cite{Os1}. Comme $\Sigma$ est de  
courbure totale finie, $\Sigma$ est conform\'ement \'equivalent au
plan complexe $\C$ (Th\'eor\`eme $9.1$ dans \cite{Os}). De plus, on peut
supposer que la sym\'etrie par 
rapport \`a $L$ corresponde sur $\C$ \`a l'application $\zeta\mapsto
\bar{\zeta}$. Sur $\C$, l'application de Gauss $g$ est d\'efinie;
$g$ est une application m\'eromorphe. Comme la normale le long de $L$
est horizontale, on a $|g(\zeta)|=1$ pour $\zeta\in\R$. De plus la
sym\'etrie par rapport \`a $L$ nous dit que $\dis g(\bar{\zeta})=
\frac{1}{\overline{g(\zeta)}}$. Comme la courbure totale est finie on  
sait enfin que $g$ est une fraction rationnelle, en regroupant tout
ces arguments on a:
$$g(\zeta)=e^{i\beta}\frac{\dis\prod_{j=1}^n(\zeta-\zeta_j)}
{\dis\prod_{j=1}^n(\zeta-\bar{\zeta_j})}$$  
avec $\beta\in[-3\pi/2,\pi/2[$.

Supposons tout d'abord que $\alpha=\pi$ et $\beta=-3\pi/2$. On
consid\`ere la suite de points du graphe de $u$ de coordonn\'ees
$(0,0,n)$. On sait que la normale au graphe en ces points converges
vers $(0,-1,0)$ lorsque $n$ tend vers $+\infty$ d'apr\`es les valeurs de $u$
sur le bord du domaine. Or si $\zeta_n\in\C$
est le point correspondant \`a $(0,0,n)$, $\zeta_n$ diverge dans $\C$
et donc $g(\zeta_n)\rightarrow e^{i\beta}=i$ ce qui nous 
donne que la valeur limite de la normale est $(0,1,0)$; on a donc une
contradiction et on suppose, dans la suite, que $\alpha\neq \pi$ ou
$\beta\neq -3\pi/2$.

On d\'efinit maintenant, sur $D(\alpha)$, la suite $(u_n)$
des dilat\'es de $u$ par $u_n(r,\theta)=\frac{1}{n}u(nr,\theta)$. On
note $\Psi_n$ la fonction conjugu\'ee de $u_n$ normalis\'ee par
$\Psi_n(0)=0$; on alors $\Psi_n\le 0$. On s'int\'eresse \`a la limite
de la normale au graphe de $u_n$ au dessus du point
$(1,\beta+\pi/2)\in D(\alpha)$ (l'appartenance \`a $D(\alpha)$ est due
\`a l'hypoth\`ese sur $\alpha$ et $\beta$). Cette normale est \'egale \`a la
normale au graphe de $u$ en $P_n$ le point au dessus de
$p_n=(n,\beta+\pi/2)$. $P_n$ diverge dans $\Sigma$ donc si $\zeta_n$
est le point de $\C$ correspondant, on a $g(\zeta_n)\rightarrow
e^{i\beta}$. Donc la normale au graphe de $u_n$ au dessus de
$(1,\beta+\pi/2)$ converge vers $(\cos\beta,\sin\beta,0)$. On a
donc une ligne de divergence qui est $L(\beta+\pi/2)$ et
$\Psi_n(1,\beta+\pi/2)\rightarrow 1$ ce qui est impossible puisque
$\Psi_n\le 0$. On vient donc de prouver que $u$ n'existe pas.   
\end{proof}

\subsection{Le cas g\'en\'eral}

On a maintenant la g\'en\'eralisation suivante du Th\'eor\`eme
\ref{T1}. 
 
\begin{thm}
Soit $\Ome$ un multi-domaine tel qu'il existe une partie $K\subset
\Ome$ et 
$\alpha\ge\pi$ tels que $\Ome\backslash K$ est isom\'etrique \`a
$D(\alpha)\backslash D(0,r)$ o\`u $D(0,r)$ est l'ensemble des point
de $D(\alpha)$ \`a distance inf\'erieure \`a $r$ de $0$. Alors il
n'existe pas de solution $u$ de l'\'equation des surfaces minimales
telle que $u$ prenne la valeur $+\infty$ (resp. $-\infty$) sur la
demi-droite correspondant \`a $L(\alpha)$ (resp. $L(-\alpha)$)
%, le
%graphe de $u$ soit de courbure totale finie
et $\Psi_u\le 0$ sur
$\Ome\backslash K$ pour un choix de normalisation pour $\Psi_u$. 
\end{thm}

\begin{proof}
La preuve est identique \`a celle du Th\'eor\`eme \ref{T1}; les seules
modifications apparaissent dans la preuve de l'\'equivalent du Lemme
\ref{L2}. On va donc r\'e\'ecrire la preuve du Lemme \ref{L2}.

On suppose \`a nouveau que $r=1$. Si le lemme n'est pas vrai on a une
suite $p_n=(r_n,\theta_n)\in D(\alpha)$ avec $r_n\rightarrow +\infty$,
$\theta_n \rightarrow \theta_\infty\in]-\alpha,\alpha[$ et
$\der{u}{\theta}(p_n) <0$. On d\'efinit alors sur $D(-\alpha+
\theta_\infty-\theta_n, \alpha+\theta_\infty-\theta_n) \backslash
D(0,{r_n}^{-1})$ la fonction $v_n$ par $v_n(r,\theta)=
{r_n}^{-1}u(r_nr,\theta-\theta_n+\theta_\infty)$. $v_n$ est une solution
de \eqref{MSE} qui vaut $+\infty$ sur $L(\alpha+\theta_\infty-\theta_n)$ et
$-\infty$ sur $L(-\alpha+ \theta_\infty-\theta_n)$. On a alors
$\der{v_n}{\theta}(1,\theta_\infty)<0$. Si on pose $\Psi_n=\Psi_{v_n}$
normalis\'ee par $\Psi_n({r_n}^{-1},\theta_\infty)={r_n}^{-1}
\Psi_u(1,\theta_n)\rightarrow 0$, on a $\Psi_n\le 0$.

On \'etudie alors les lignes de divergence de la suite $(v_n)$ sur le
multi-domaine limite qui est $D(\alpha)$. Comme plus haut, les seules
lignes de divergence sont les $L(\gamma)$ ou les droites incluses dans
$D(\alpha)$. Si $L(\gamma)$ et une ligne de divergence on sait que pour
une sous suite $\Psi_{n'}\rightarrow \Psi$ avec $\Psi(0)=0$ et
$\Psi(\rho,\gamma)= \pm \rho$. Or $\Psi\le 0$ donc
$\Psi(\rho,\gamma)=-\rho$ et la normale limite le long de $L(\gamma)$
est $(-\sin\gamma, \cos\gamma)$. Il y a donc une seule normale limite
possible. Comme dans la preuve originale ceci implique que les
$L(\gamma)$ sont les seules lignes de divergence possibles.

Comme $\der{v_n}{\theta}(1,\theta_\infty)\le 0$, le point
$(1,\theta_\infty)$ appartient \`a $\boB(v_n)$ le domaine de
convergence de $v_n$. La composante connexe de $\boB(v_n)$ contenant
$(1,\theta_\infty)$ est alors $D(\gamma_1,\gamma_2)$ pour
$\gamma_1<\gamma_2$. Soit alors $v$ une limite d'une sous-suite de $(v_n)$
sur $D(\gamma_1,\gamma_2)$,
on a $\Psi_v=\lim \Psi_{n'} \le 0$ et $v$ vaut $+\infty$ sur
$L(\gamma_2)$ et $-\infty$ sur $L(\gamma_1)$.
% Maintenant on sait que la
%courbure totale du graphe de $v_n$ ne d\'epend pas de $n$ et vaut la
%courbure totale du graphe de $u$ au dessus de $\Ome \backslash
%K$. Ainsi par le lemme de Fatou, la courbure totale du graphe de $v$
%est inf\'erieure \`a la limite inf\'erieure des courbures totales des
%graphes de $v_{n'}$ au dessus de $D(\gamma_1,\gamma_2)$ qui est
%elle-m\^eme inf\'erieure \`a la courbure totale du graphe de $u$.
Ainsi  $v$
est une solution de \eqref{MSE} qui contredit la Proposition
\ref{P4}. L'\'equivalent du Lemme \ref{L2} est donc prouv\'e.   
\end{proof}

%%% Local Variables: 
%%% mode: latex
%%% TeX-master: t
%%% End: 